\documentclass{scrartcl}
\usepackage[utf8]{inputenc}
\usepackage{lineno,hyperref}
\usepackage{amsmath}
\usepackage{amssymb}
\usepackage{color}
\usepackage{comment}
\usepackage{tikz}
\usepackage{pgfplots}
\pgfplotsset{compat=1.15}
\usepackage{float,todonotes}
\usetikzlibrary{patterns}
\modulolinenumbers[5]
\usepackage{enumerate}

 \newtheorem{theorem}{Theorem}[section]

  \newtheorem{definition}[theorem]{Definition}
  \newtheorem{remark}[theorem]{Remark}

\usepackage{bbm}
\newcommand{\real}[0]{\mathbbm{R}}
\newcommand{\nat}[0]{\mathbbm{N}}
\usepackage[normalem]{ulem}
\usepackage{cancel}
\newcommand{\yS}[0]{y^{\{S\}}}

\newcommand{\yShat}[0]{\hat{y}^{\{S\}}}
\newcommand{\ySbar}[0]{\bar{y}^{\{S\}}}
\newcommand{\yF}[0]{y^{\{F\}}}

\newcommand{\yFtilde}[0]{\tilde{y}^{\{F\}}}
\newcommand{\yFbar}[0]{\bar{y}^{\{F\}}}
\newcommand{\Oh}[0]{\mathcal{O}}
\definecolor{darkgreen}{rgb}{0.05,0.5,0.05}

\begin{document}

\title{Spline-oriented inter/extrapolation-based multirate schemes of higher order}

\author{Kevin Schäfers, Andreas Bartel, Michael Günther \\
	and Christoph Hachtel \\[2ex]
	\small 
University of Wuppertal, \\ 
	\small 
Faculty of Mathematics and natural sciences,  IMACM,\\ 
	\small 
40097 Wuppertal (Germany) \\
	\small 
email: \{schaefers,bartel,guenther,hachtel\}@math.uni-wuppertal.de}
\date{}
 \maketitle

\begin{abstract}
	\textbf{Abstract.} 
Multirate integration uses different time step sizes for different 
components of the solution based on the respective transient behavior.
For inter/extrapolation-based multirate schemes, we
construct a new subclass of schemes by using clamped cubic splines to obtain multirate schemes up to order 4. 
Numerical results 
for a $n$-mass-oscillator demonstrate that 4th order of convergence can be achieved for this class of schemes.
\end{abstract}


\medskip

{\small \noindent\textbf{Acknowledgements: } The third author indebted to the funding provided by the EU Horizon 2020 Marie Skłodowska-Curie grant agreement no. 765374.}

\section{Introduction}
\noindent Many technical applications, e.g. electric circuits~\cite{gear1984multirate,Bartel2022},
can be modeled as coupled systems of ordinary differential equations (ODEs). 
Often, the transient behavior of these applications is characterized by different time constants. %
At a given time point $t_n$, there are many slowly evolving components $y^{\{S\}} \in \real^{d_S}$ and a few 
fast components $y^{\{F\}} \in \real^{d_F}$ ($d_S \gg d_F$). To exploit this,
multirate integration schemes are developed following the pioneer work~\cite{rice1960split}. 
Thereby, the fast part $y^{\{F\}}$ 
is integrated with a small step size $h$ (micro step) and the slow part $y^{\{S\}}$ 
using a larger step size $H = m\cdot h$ ($m \in \nat$), (macro step). The appropriate choice of the coupling variables is the challenging part. Here, we focus on inter/extrapolation-based multirate schemes to describe the novel construction of spline-oriented multirate schemes of higher order.
Numerical results demonstrate an order four convergence based on cubic spline coupling.

\section{Inter/extrapolation-based multirate schemes}
\label{ch:MultirateSchemes}
\noindent We consider the component-wise partitioned initial value problem (IVP)
\begin{align}
    \begin{aligned}\label{IVP}
    \dot{y}^{\{S\}} &= f^{\{S\}}\left(t,y^{\{S\}},y^{\{F\}}\right), 
      \quad & y^{\{S\}}(t_0) &= y_{0}^{\{S\}} \in \real^{d_S}, 
      \\
    \dot{y}^{\{F\}} &= f^{\{F\}}\left(t,y^{\{S\}},y^{\{F\}}\right), 
      \quad
    & y^{\{F\}}(t_0) &= y_{0}^{\{F\}} \in \real^{d_F},
\end{aligned}
\end{align}
of coupled ODEs on $[t_0,T]$, where $f^{\{S\}}$ and $f^{\{F\}}$ are 
sufficiently smooth.%

Depending on the sequence of computation, one can distinguish five versions of inter/ extrapolation-based multirate schemes: \textit{fully-decoupled approach} \cite{Bartel2022}, \textit{decoupled slowest-first approach} \cite{gear1984multirate}, \textit{decoupled fastest-first approach} \cite{gear1984multirate}, \textit{coupled first-step approach} \cite{gunther2001multirate} and the \textit{coupled slowest-first approach}, which is a special case of the time-stepping strategy introduced in \cite{savcenco2007multirate} and traces back to Rice \cite{rice1960split}. In this work, we focus on the decoupled slowest-first (DSF) approach.

\definition[DSF approach]%
{%
We advance the solution of~\eqref{IVP} from $t_n$ to $t_n+H$.
Firstly, the slow subsystem is numerically integrated 
by a scheme $\Phi_H$ and an extrapolated waveform $\tilde{y}^{\{F\}}$ for 
$y^{\{F\}}$.
a scheme $\hat{\Phi}_h$ and an interpolated waveform $\hat{y}^{\{S\}}$ for the slow components based on information from the current time window $[t_n,t_n+H]$.
}\normalfont

\remark{All five inter/extrapolation-based multi\-rate schemes have convergence order $p$, if the basic integration schemes 
have convergence order $p$ and the inter/extrapolation schemes are of approximation order $p-1$ \cite{Bartel2020}.}
\normalfont

\section{Multirate schemes using spline-oriented inter/extrapolation}
\label{ch:Splines}
\noindent As derivative information is provided by the ODE, we consider clamped cubic splines of order three~\cite{bulirsch2002introduction}.
This enables multirate schemes of order $p=4$. 

\paragraph{DSF clamped-spline approach} We have two parts:
\begin{enumerate}
    \item[(i)] 
    After a numerical approximation is obtained for a macro step
    $[t_{n-1},t_{n}]$, 
    compute a clamped cubic spline $S^{\{F\}}_{n}(t)$ for the fast components 
    based on computed micro step 
    values $y_{n-1+i/m}^{\{F\}}$, $i=0,\ldots,m$,
    and the 
    derivative information
    \[
       \dot{y}_{n-1}^{\{F\}} = f^{\{F\}}\left(t_{n-1},y_{n-1}^{\{S\}},y_{n-1}^{\{F\}}\right),
       \qquad 
       \dot{y}_{n}^{\{F\}} = f^{\{F\}} \left(t_{n},y_{n}^{\{S\}},y_{n}^{\{F\}}\right).
    \]   
    \item[(ii)] After the integration of the slow subsystem, compute a cubic polynomial $S^{\{S\}}_{n+1}(t)$ for the slow components by using the idea of clamped cubic splines based on the computed values $y_n^{\{S\}}$ and $y_{n+1}^{\{S\}}$, 
    and the derivative information 
    $$
      \dot{y}_n^{\{S\}} = f^{\{S\}}\left(t_n,y_n^{\{S\}},y_n^{\{F\}}\right) \quad\text{and}\quad  \dot{y}_{n+1}^{\{S\}} = f^{\{S\}}\left(t_{n+1},y_{n+1}^{\{S\}},\tilde{y}_{n+1}^{\{F\}}\right).
    $$
    Note that we have to use the extrapolated value $\tilde{y}_{n+1}^{\{F\}}$ as $y_{n+1}^{\{F\}}$ is unknown. Fortunately (except of some special cases \cite{gear1984multirate}), the multirate 
    setting implies that $\left\lVert \frac{\partial f^{\{S\}}}{\partial y^{\{F\}}}\right\rVert$ is small such that the use of the extrapolated value is not problematic.
\end{enumerate}

\paragraph{The first macro step}
We need a special treatment for the first macro step as there exists no spline $S_0^{\{F\}}(t)$. Thus, for the first macro step $[t_0,t_0 + H]$, we apply $m$ micro steps of a singlerate integration scheme. At the end of this macro step, one has to compute the clamped cubic spline $S_1^{\{F\}}(t)$ on $[t_0,t_1]$.

\paragraph{The general procedure} Con\-si\-der the time window $[t_n,t_{n+1}], \, n \geq 1$.
\begin{enumerate}[(i)]%
    \item integrate the slow subsystem for $y^{\{S\}}$ with step size $H$ using the extra\-po\-lated waveform $\left.\tilde{y}^{\{F\}} = S^{\{F\}}_{n}\right\vert_{[t_{n}-h,t_{n}]}$,
    \item compute the cubic polynomial $S_{n}^{\{S\}}(t)$ on $[t_n,\,t_{n+1}]$,
    \item perform the $m$ micro steps for the fast subsystem for $y^{\{F\}}$ with step size $h$ using the interpolated waveform $\hat{y}^{\{S\}} = S_{n}^{\{S\}}$,
    \item compute the clamped cubic spline $S_{n+1}^{\{F\}}(t)$ on $[t_{n},\, t_{n+1}]$ with nodes 
    \mbox{$t_n\!+\!ih$}, $i=0,\ldots,m$.
\end{enumerate}

\noindent These methods are referred to as \emph{spline-oriented multirate schemes}.
\paragraph{Convergence Analysis}
Let the solution be computed on $[t_0,t_n]$. Then, one computes a cubic polynomial  $\tilde{y}^{\{F\}}$  based on the fast data
    \[ 
      y_{n-1+i/m}^{\{F\}},\; 
      i=0,\ldots,m,\;
      f^{\{F\}}\left(t_{n-1},\; y_{n-1}^{\{S\}},y_{n-1}^{\{F\}}\right),\; 
      f^{\{F\}}\left(t_n, y_{n}^{\{S\}},y_n^{\{F\}}\right)
    \]
    of the last macro step $[t_{n-1},t_n]$.    
    Thus it holds
      $\lvert \yF(t) - \yFtilde(t)\rvert = \Oh(H^4)$ for $t \in [t_n,t_{n+1}]$.
    For the computation of the next macro step $[t_n, t_{n+1}]$, the use of the extrapolated waveform in step (i) of the general procedure results in a perturbed, decoupled ODE system for the slow subsystem. This introduces a model error to the slow part
    \begin{align}\label{modelErrSlow}
         \dot{\bar{y}}^{\{S\}} ={f}^{\{S\}}\left(t,\bar{y}^{\{S\}},\tilde{y}^{\{F\}} (t)\right)
         \;= \; {f}^{\{S\}}\left(t,\bar{y}^{\{S\}}, {y}^{\{F\}} \right) + 
           \Psi_n^{\{S\}}(t) 
    \end{align}
    with $\Psi_n^{\{S\}} \in \Oh(H^4)$.
    Subtracting this from the original ODE for $\yS$ and using Gronwall's lemma, we deduce
    \[
      \left\vert \yS(t) - \ySbar(t) \right\vert = \Oh(H^4), \quad t\in [t_n,t_{n+1}]    .
    \]
    An integration scheme of order $p \ge 4$ yields an approximation $y_{n+1}^{\{S\}}$
    with error
    \begin{equation}\label{eq:modelling_error}
     \left\vert\yS(t_{n+1}) - y_{n+1}^{\{S\}}\right\vert \le 
     \left\vert \yS(t_{n+1}) - \ySbar(t_{n+1}) \right\vert + \left\vert\ySbar(t_{n+1}) - y_{n+1}^{\{S\}}\right\vert, 
    \end{equation}
    i.e., the error is, like the approximation and model error, in $\Oh(H^4)$.
    To update the fast part, step (ii) computes a cubic polynomial $\yShat$ using  
     \[
        y_{n}^{\{S\}}, \,
        y_{n+1}^{\{S\}},\, 
        f^{\{S\}}\left(t_n,y_{n}^{\{S\}},y_{n}^{\{F\}}\right),\, f^{\{S\}}\left(t_{n+1},y_{n+1}^{\{S\}},\tilde{y}_{n+1}^{\{F\}}\right).
     \]
     As 
     $\tilde{y}_{n+1}^{\{F\}} = y_{n+1}^{\{F\}} + \Oh(H^4)$,
     it holds $\lvert y^{\{S\}}(t) - \hat{y}^{\{S\}}(t)\rvert = \Oh(H^4)$ for $t \in [t_n,t_{n+1}]$.
     This leads to the perturbed, decoupled fast ODE system
     ($\Psi_n^{\{F\}}(t)  \in \Oh(H^4)$)
    \begin{align}\label{modelErrFast}
        \dot{\bar{y}}^{\{F\}} =    {f}^{\{F\}}\left(t,\hat{y}^{\{S\}},\bar{y}^{\{F\}}\right)
        =
          {f}^{\{S\}}\left(t, {y}^{\{S\}},\bar{y}^{\{F\}}\right)
          +  \Psi_n^{\{F\}}(t)
    \end{align}
    In step (iii), computing $m$ micro steps for \eqref{modelErrFast} with a numerical integration scheme of order $p \geq 4$ yields a fourth order update as it holds 
    \begin{equation}\label{eq:modelling_errorFast}
        \left\vert\yF(t_{n+1}) - y_{n+1}^{\{F\}}\right\vert \le 
     \left\vert \yF(t_{n+1}) - \yFbar(t_{n+1}) \right\vert + \left\vert\ySbar(t_{n+1}) - y_{n+1}^{\{F\}}\right\vert.
    \end{equation}
In the end, the spline-oriented approximation at the final time $t=T$ can be viewed as a numerical approximation 
of the decoupled ODE system, where the coupling variables are represented by inter- and extrapolated functions of $t$. This yields modifications of the right hand sides: $\Psi^{\{S\}}$ and $\Psi^{\{F\}}$. Both functions are successively being build up during the integration and result in a modelling error of order 4.
This gives: 

\theorem[Convergence]{The spline-oriented multirate scheme is a convergent procedure. The order depends on the employed basic schemes and spline approximation. In the case of clamped cubic splines and a numerical integration scheme of order $p \geq 4$, we obtain a method of order four.
}
\remark{The spline-oriented multirate scheme has a predictor-corrector type of structure.}
 
\definition[Spline-oriented multirate RK scheme]
{\label{exampleScheme}
We consider a 
$s$-staged Runge--Kutta (RK) scheme~\cite{bulirsch2002introduction} given by the Butcher tableau $\left(A,b,c\right)$. Based on the DSF-approach, one macro step of a spline-oriented multirate RK scheme with $m$ micro steps, applied to \eqref{IVP}, advances the solution $\left(y_{n}^{\{S\}},y_{n}^{\{F\}}\right)$ at $t_n = t_0 + n \cdot H$ to the solution $\left(y_{n+1}^{\{S\}},y_{n+1}^{\{F\}}\right)$ at $t_{n+1} = t_{n} + H$ as follows:\\
1.) macro step
\begin{subequations}\label{MR-RKV}
\begin{align}
\begin{split}\label{macro-step}
    y_{n+1}^{\{S\}} &= y_{n}^{\{S\}} + H \sum\limits_{i=1}^s b_i k_{i}^{\{S\}}, \\
    k_{i}^{\{S\}} &= f^{\{S\}}\Bigl( t_n + c_i H,\,\, y_{n}^{\{S\}} + H \sum \limits_{j=1}^s a_{i,j} k_{j}^{\{S\}},\tilde{y}_{n+c_i}^{\{F\}}\Bigr), \\
\end{split}
\intertext{2.) $m$ micro steps (for $i=1,\ldots,s$ and $\lambda = 0,\ldots,m-1$)}
\begin{split}\label{micro-steps}
    y_{n+\frac{\lambda+1}{m}}^{\{F\}} &= y_{n+\frac{\lambda}{m}}^{\{F\}} + h \sum\limits_{i=1}^s b_i k_{i}^{\{F,\lambda\}}, \\
    k_i^{\{F,\lambda\}} &= f^{\{F\}}\Bigl(t_n + (\lambda + c_i)h, \,\, \hat{y}_{n+\frac{\lambda+c_i}{m}}^{\{S\}},y_{n+\frac{\lambda}{m}}^{\{F\}} + h \sum\limits_{j=1}^s a_{i,j} k_j^{\{F,\lambda\}}\Bigr).
\end{split}
\end{align}
\end{subequations}
}\normalfont

\section{Numerical results}
\label{ch:Results}


%
%
%
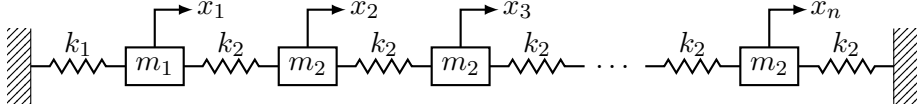
\begin{figure}[bth]
    \centering   
    \begin{tikzpicture} 
  \def\springlength{1.25} 
  \tikzset{ 
      box/.style={draw,outer sep=0pt,thick}, 
      spring/.style={thick,decorate, 
        decoration={zigzag,pre length=0.25cm,post length=0.25cm,segment length=6}}, 
      springdescr/.style={yshift=.3cm}, 
      directiondescr/.style={xshift=.75cm},
      forcedescr/.style={yshift=-.3cm},
      ground/.style={box,draw=none,fill, 
        pattern=north east lines,minimum width=0.3cm,minimum height=1cm}, 
      mass/.style={box,minimum width=.5cm,minimum height=.5cm}, 
      direction/.style={-latex,thick} 
    } 
  \node (wall) [ground] {}; 
  \draw (wall.south east) -- (wall.north east); 
 \draw [spring] (wall.east) -- node[springdescr]{$k_1$} 
        +(\springlength,0)node[mass,anchor=west](M1){$m_1$}; 
 \draw [spring] (M1.east) -- node[springdescr]{$k_2$} 
        +(\springlength,0)node[mass,anchor=west](M2){$m_2$}; 
 \draw [spring] (M2.east) -- node[springdescr]{$k_2$} 
        +(\springlength,0)node[mass,anchor=west](M3){$m_2$}; 
 \draw [spring] (M3.east) -- node[springdescr]{$k_2$} 
        +(\springlength,0)node[anchor=west](M4){$\dotsc$}; 
 \draw [spring] (M4.east) -- node[springdescr]{$k_2$} 
        +(\springlength,0)node[mass,anchor=west](Mn){$m_2$}; 
 		
 \draw [direction] (M1.north) |- node[directiondescr]{$x_1$} +(0.5,0.5);
  \draw [direction] (M2.north) |- node[directiondescr]{$x_2$} +(0.5,0.5);
 \draw [direction] (M3.north) |- node[directiondescr]{$x_3$} +(0.5,0.5);
  \draw [direction] (Mn.north) |- node[directiondescr]{$x_n$} +(0.5,0.5);
 \begin{scope}[xshift=11.65cm] 
 \node (wall) [ground] {}; 
 \draw (wall.south west) -- (wall.north west); 
  \draw [spring] (wall.west) -- node[springdescr]{$k_2$} 
        +(-\springlength,0); 
 \end{scope} 
\end{tikzpicture}
\caption{System with $n$ masses: one light mass $m_1$, $n\!-\!1$ heavy masses; $n\!+\!1$ springs, one strong spring $k_1$, $n$ light springs $k_2$. The system is attached to fixed walls.}
\label{fig:massSpring_multi}
\end{figure}
\begin{figure}[H]
\begin{minipage}{0.49\textwidth}
    \centering
    \resizebox{\textwidth}{!}{%
    \input{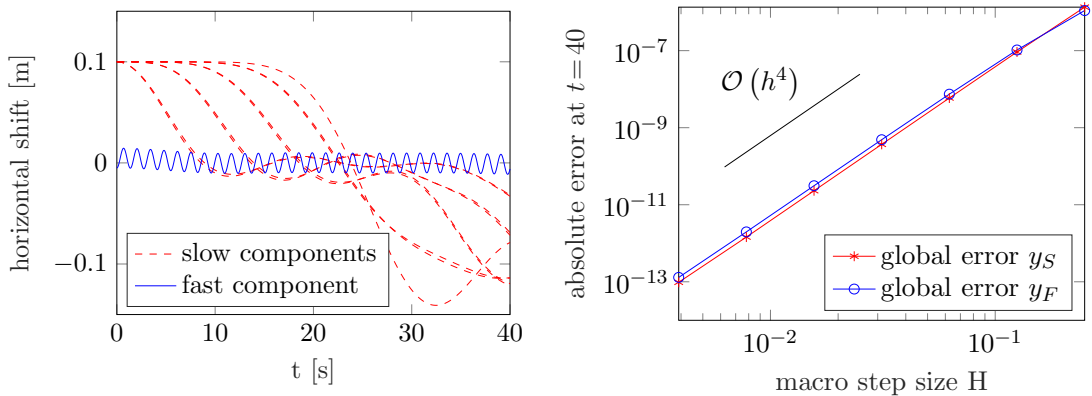}%
    }
\end{minipage}
\begin{minipage}{0.49\textwidth}
    \centering
    \resizebox{\textwidth}{!}{%
%
%
\begin{tikzpicture}

\begin{axis}[%
width=2.25in,
height=1.75in,
at={(0.758in,0.481in)},
scale only axis,
xmode=log,
xmin=0.00390625,
xmax=0.25,
xminorticks=true,
xlabel style={font=\color{white!15!black}},
xlabel={macro step size H},
ymode=log,
ymin=1e-14,
ymax=1.34530762800305e-06,
yminorticks=true,
ylabel style={font=\color{white!15!black}},
ylabel={absolute error at $t\!=\!40$},
axis background/.style={fill=white},
legend style={at={(0.97,0.03)}, anchor=south east, legend cell align=left, align=left, draw=white!15!black}
]
\addplot [color=red, mark=asterisk, mark options={solid, red}]
  table[row sep=crcr]{%
0.25	1.34530762800305e-06\\
0.125	9.22126026047759e-08\\
0.0625	5.88388260305742e-09\\
0.03125	3.69171581155157e-10\\
0.015625	2.30839094883698e-11\\
0.0078125	1.44704482570301e-12\\
0.00390625	9.7017024975559e-14\\
0	0\\
0	0\\
0	0\\
};
\addlegendentry{$\text{global error } y_S$}

\addplot [color=blue, mark=o, mark options={solid, blue}]
  table[row sep=crcr]{%
0.25	1.09767223397474e-06\\
0.125	1.03995264344255e-07\\
0.0625	7.41079278553534e-09\\
0.03125	4.8719849220833e-10\\
0.015625	3.11263411805029e-11\\
0.0078125	1.97134172995121e-12\\
0.00390625	1.29875047580304e-13\\
0	0\\
0	0\\
0	0\\
};
\addlegendentry{$\text{global error } y_F$}

\addplot [color=black, forget plot]
  table[row sep=crcr]{%
0.025	2.44140625e-08\\
0.0125	1.52587890625e-09\\
0.00625	9.5367431640625e-11\\
};
\node[above left] at (axis cs: 0.015,3e-9) {$\mathcal{O}\left(h^4\right)$};
\end{axis}

\begin{axis}[%
width=2.9in,
height=2.15in,
at={(0in,0in)},
scale only axis,
xmin=0,
xmax=1,
ymin=0,
ymax=1,
axis line style={draw=none},
ticks=none,
axis x line*=bottom,
axis y line*=left,
legend style={legend cell align=left, align=left, draw=white!15!black}
]
\end{axis}
\end{tikzpicture}%
%
    }
\end{minipage}
\caption{Numerical results for a  non-stiff  10-mass-oscillator~\eqref{eq:nMassSpringODE} for MR-RK4. Left: numerical solutions of slow (red) and fast (blue) components. Right: convergence order for fixed multirate factor $m\!=\!20$ ($H=mh$) as absolute error at $t\!=\!40$. }
\label{fig:Results}
\end{figure}

 
\noindent We consider the line configuration of $n$ masses and $n\!+\!1$ springs as given in Fig.~\ref{fig:massSpring_multi}. The equation of motion reads ($x_i$ position of $i$th mass):
\begin{align}
\label{eq:nMassSpringODE}
\begin{split}
m_1 \ddot{x}_1 &= -(k_1+k_2)x_1 + k_2 x_2, \\
m_2 \ddot{x}_2 &= k_2 x_1 - 2k_2 x_2 + k_2 x_3, \\
&\phantom{I}\vdots \\
m_2 \ddot{x}_{n-1} &= k_2 x_{n-2} - 2k_2 x_{n-1} + k_2 x_n, \\
m_2 \ddot{x}_n &= k_2 x_{n-1} - 2k_2 x_n.
\end{split}
\end{align}
\noindent To demonstrate the construction of an order four multirate scheme, we use as underlining schemes~\ref{exampleScheme} the classical Runge-Kutta (RK4) scheme~\cite{bulirsch2002introduction} which we refer to as
MR-RK4.
With multirate factor $m=20$, we solve the ODE system \eqref{eq:nMassSpringODE} for $n=10$ numerically on the time interval $[0,40]$ with masses $\left(m_1,m_2\right) = \left(1,20\right)$, spring constants $\left(k_1,k_2\right) = \left(20,1\right)$ and initial values 
\begin{align*}
    \begin{pmatrix}
    x_1(0),\;
    \dot{x}_1(0)
    \end{pmatrix} = \begin{pmatrix}
    -0.005,\;
    0
    \end{pmatrix}, \quad \begin{pmatrix}
    x_i(0),\; 
    \dot{x}_i(0)
    \end{pmatrix} = \begin{pmatrix}
    0.1,\;
    0
    \end{pmatrix}, \quad \forall i=2,\ldots,n.
\end{align*}
Numerical results are shown in Figure \ref{fig:Results}. Order 4 
of the MR-RK4 is observed.
\noindent We discuss the achieved accuracy of MR-RK4 ($m=20$) compared to the singlerate case (RK4 with micro step size $h$) in the above test~\eqref{eq:nMassSpringODE}: 
\noindent i)     The error of the fast part $y_F$ is above three orders of magnitude less than the absolute error of the singlerate RK4 scheme. 
ii) The error of the slow part $y_S$ is above 1.5 orders of magnitude less than the 
    absolute error of the singlerate RK4 scheme.
On the other hand, multirate saves right-hand evaluations. For the ODE system \eqref{eq:nMassSpringODE}, for one macro step, the number of scalar function evaluations is:
\[
\text{singlerate: }  4 \cdot m \cdot n \quad \text{vs.}\quad 
\text{multirate: }  4 \cdot \left(m \cdot 1 + 1 \cdot (n-1)\right).
\]
For $m=20$ and $n=10$, we thus have 800 vs. 116 evaluations.

\section{Conclusion and outlook}
\noindent We combined successfully inter/extrapolation-based multirate schemes with a clamped cubic spline. Thereby, we preserve the order of the basic integration scheme.
Next steps are some further stability investigations 
and generalizations of this class of multirate schemes. 
Apart from this, we then plan 
to extend this strategy to the application of DAEs. Of course, the interpolation of the algebraic variables will move us from the manifold. This needs to be treated.

\bibliography{mybibfile}

\end{document}